\documentclass[12pt]{amsart}
\usepackage{enumerate}
\usepackage{amsmath}
\usepackage{amssymb}
\usepackage{amsfonts}
\usepackage{amsthm, upref}
\usepackage{graphicx}
\usepackage[usenames, dvipsnames]{color}

    \numberwithin{equation}{section}

\usepackage{bm}
         
         \bmdefine\alphab{\mathbf{\alpha}}
\bmdefine\betab{\mathbf{\beta}}
\bmdefine\sigmab{\mathbf{\sigma}}

\newcommand{\comment}[1]{}
\newcommand{\eq}{\begin{equation}}
\newcommand{\en}{\end{equation}}
\newcommand{\pp}{\mathbb{P}}
\newcommand{\rr}{\mathbb{R}}
\newcommand{\nn}{\mathbb{N}}

\newcommand{\ep}{\hfill $\Box$}

\begin{document}

\theoremstyle{plain}
\newtheorem{thm}{Theorem}
\newtheorem{lemma}[thm]{Lemma}
\newtheorem{prop}[thm]{Proposition}
\newtheorem{cor}[thm]{Corollary}

\theoremstyle{definition}
\newtheorem{defn}{Definition}
\newtheorem{cond}{Condition}
\newtheorem{asmp}{Assumption}
\newtheorem{notn}{Notation}
\newtheorem{prb}{Problem}

\theoremstyle{remark}
\newtheorem{rmk}{Remark}
\newtheorem{exm}{Example}
\newtheorem{clm}{Claim}

\title[Markov chains]{Some universal estimates for reversible Markov chains}

\author{Mykhaylo Shkolnikov}
\address{Mathematical Sciences Research Institute \\ 17 Gauss Way \\ Berkeley, CA 94720-5070}
\email{mshkolni@gmail.com}

\thanks{This research was partially supported by NSF grant DMS-08-06211}

\date{\today}

\begin{abstract}
We obtain universal estimates on the convergence to equilibrium and the times of coupling for continuous time irreducible reversible finite-state Markov chains, both in the total variation and in the $L^2$ norms. The estimates in total variation norm are obtained using a novel identity relating the convergence to equilibrium of a reversible Markov chain to the increase in the entropy of its one-dimensional distributions. In addition, we propose a universal way of defining the ultrametric partition structure on the state space of such Markov chains. Finally, for chains reversible with respect to the uniform measure, we show how the global convergence to equilibrium can be controlled using the entropy accumulated by the chain.    
\end{abstract}

\maketitle

\section{Introduction}

Recently, the convergence to equilibrium of slowly mixing Markov chains appearing in statistical physics has attracted much attention. In this framework continuous time irreducible reversible Markov chains are defined by choosing the transition rates from a state (usually, a spin configuration) $a$ to a state (spin configuration) $b$ to be proportional to $e^{-\beta(E(b)-E(a))_+}$, where $E$ is an energy functional and $\beta$ stands for the inverse temperature, which in this context is chosen to be large: $\beta\gg1$. In this low temperature regime, a recurring feature is that the energy landscape given by $E$ divides the state space into sets of metastable (or, stable) states, which are separated by potential wells. The convergence to equilibrium of the corresponding Markov chain, started in a metastable state, is then governed by the time it takes to overcome the respective potential wells in order to reach the part of the state space with the lowest energy.

\bigskip

The potential theoretic approach to metastability developed in the articles \cite{BEGK1}, \cite{BEGK2} and \cite{BEGK3} (see also the excellent summaries \cite{B1} and \cite{B2}) has been used to obtain precise information on metastable transitions for reversible Markov chains associated with several models of statistical physics. These include certain disordered mean field models (see \cite{BEGK2}) and, more specifically, the Curie-Weiss model with a random field taking finitely many values (see \cite{BEGK2}, \cite{BBI1} and \cite{BBI2}). Other examples of slowly mixing reversible Markov chains, in which the metastable behavior has been analyzed in detail, include the Glauber dynamics for the two-dimensional Ising model on a torus and its generalizations (see \cite{NS1} and \cite{NS2}), the three-dimensional Ising model on a torus (see \cite{BC}) and the classical Curie-Weiss model (see \cite{CGOV} and \cite{LLP}). Moreover, the first exit problem from a domain for reversible chains with exponentially small transition probabilities was studied in the article \cite{OS}. In a different but related line of research, initiated by the article \cite{K}, the metastable transitions are studied for diffusions with a small diffusion parameter, which are confined in a potential having several local minima (see \cite{BEGK4} and \cite{BGK} for a recent account on this problem). 

\bigskip

Here, we take a different viewpoint. Instead of analyzing a specific Markov chain in detail, we try to understand some universal aspects of the ultrametric structure, that is, the presence of multiple time scales in a general irreducible reversible finite-state Markov chain. We obtain universal estimates on the convergence to equilibrium and the times of coupling in this abstract framework. We prove such results both in the total variation norm and in the $L^2$ norm. In the case of the total variation norm, we utilize a novel entropy identity relating the convergence to equilibrium of the chain to the increase of the entropy of its one-dimensional distributions. 

\bigskip

In addition, we propose a universal way of defining the ultrametric partition structure, that is, a sequence of partitions of the state space corresponding to the different time scales on which convergence to equilibrium occurs. Finally, in the case that the chain is reversible with respect to the uniform measure, we show how the entropy of the one-dimensional distributions can be utilized to control the global convergence to equilibrium of the chain.  

\bigskip

To give examples of the type of results we obtain, we introduce a set of notations. Let $X$ be a continuous time irreducible reversible Markov chain on a set $I=\{1,2,\ldots,n\}$ of $n$ elements. Moreover, for an $a\in I$ and a $\,t\geq0$ let $P^a_t$ be the one-dimensional distribution of the chain at time $t$, when started in $a$. Finally, write $\|.\|_{TV}$ for the total variation norm, $\nu$ for the invariant distribution of $X$ and let $H(\nu)$ be the entropy $-\sum_{a\in I} \nu(a)\log\nu(a)$ of the invariant distribution $\nu$.  

\bigskip

Before stating our first result rigorously, we would like to provide the reader with some intuition by giving an example. Fix natural numbers $M,\,N\geq3$ and consider the graph given by an arrangement of $M$ $N$-cycles in a cycle of size $M$. Now, let $X$ be the continuous time Markov chain on this graph, which has a transition rate $\rho_1>0$ for neighboring vertices belonging to the same $N$-cycle and a transition rate $\rho_2>0$ for neighboring vertices belonging to different $N$-cycles. Since the generating matrix of this Markov chain is symmetric, it is reversible with respect to the uniform distribution on the set of vertices of the graph. Next, suppose that $\rho_2$ is much smaller than $\rho_1$. Then, it is intuitively clear that, for every vertex $a$ of the described graph, the quantity $\|P^a_{3t}-P^a_t\|_{TV}$ can be only large on the two disjoint time intervals, during which the Markov chain mixes on the $N$-cycle containing $a$ and on the $M$-cycle comprised by the $M$ $N$-cycles, respectively. Under the scale-invariant measure, which has the density $\frac{1}{t}$ on the time axis $[0,\infty)$, the union of these two time intervals has a total measure of order $\log N + \log M=\log(MN)$. Thus, it is logarithmic in the size of the state space of $X$. The purpose of Theorems \ref{thmTVgen} and \ref{thm1} below is to show that the latter property is universal for continuous time irreducible reversible Markov chains, in the sense that the order of magnitude in this example is an upper bound on the size of the corresponding quantity for a general reversible Markov chain.

\smallskip

\begin{thm}\label{thmTVgen}
Let $X$ be a continuous time irreducible reversible Markov chain on the set $I=\{1,2,\ldots,n\}$ and let $\nu$ be its invariant distribution. Then, the following is true.
\begin{enumerate}[(a)]
\item For every $\delta>0$, there is a constant $C(\delta)>0$ depending only on $\delta$ (and not on $n$ or the particular Markov chain) such that 
\eq
\sum_{a\in I}\nu(a)\,\int_0^\infty \mathbf{1}_{\{t\geq0:\;\|P^a_{3t}-P^a_t\|_{TV}\geq\delta\}}\,\frac{1}{t}\,\mathrm{d}t \leq C(\delta)\,H(\nu).
\en
In particular, for every $\delta,\,\epsilon>0$, there exists a constant $C_\epsilon(\delta)>0$ depending only on $\delta$ and $\epsilon$ (and not on $n$ or the particular Markov chain) such that
\eq
\nu\left(a\in I:\;\int_0^\infty \mathbf{1}_{\{t\geq0:\;\|P^a_{3t}-P^a_t\|_{TV}\geq\delta\}}\,\frac{1}{t}\,\mathrm{d}t \geq C_\epsilon(\delta)\,H(\nu)
\right)<\epsilon.
\en 
\item For every $\delta>0$, there is a constant $\tilde{C}(\delta)>0$ depending only on $\delta$ (and not on $n$ or the particular Markov chain) such that
\eq
\sum_{(a,b)\in I^2}\nu(a)\,\nu(b)\,\int_0^\infty \mathbf{1}_{\{t\geq0:\;\|P^a_{t}-P^b_{t}\|_{TV}-\|P^a_{3t}-P^b_{3t}\|_{TV}\geq\delta\}}\,\frac{1}{t}\,\mathrm{d}t \leq \tilde{C}(\delta)\,H(\nu).
\en
\end{enumerate}
In particular, for every $\delta,\,\epsilon>0$, there exists a constant $\tilde{C}_\epsilon(\delta)>0$ depending only on $\delta$ and $\epsilon$ (and not on $n$ or the particular Markov chain) such that
\eq
(\nu\times\nu)\left((a,b)\in I^2:\int_0^\infty \mathbf{1}_{\{t\geq0:\|P^a_{t}-P^b_{t}\|_{TV}-\|P^a_{3t}-P^b_{3t}\|_{TV}\geq\delta\}}\frac{1}{t}\mathrm{d}t\geq \tilde{C}_\epsilon(\delta)\,H(\nu)
\right)<\epsilon.
\en
\end{thm}

\smallskip

We remark at this point that universal estimates as in Theorem \ref{thmTVgen} can only be obtained under the \textit{scale-invariant measure} $\frac{1}{t}\,\mathrm{d}t$ on the time axis $[0,\infty)$, which has the property
\eq
\int_{t_1}^{t_2} \frac{1}{t}\,\mathrm{d}t=\int_{\eta t_1}^{\eta t_2} \frac{1}{t}\,\mathrm{d}t 
\en
for all $\eta>0$ and $0<t_1<t_2<\infty$. This can be easily seen by slowing down or speeding up the chain by a constant factor. 

\bigskip

To give an example of a result in the framework of $L^2$ convergence, a set of auxiliary notations is needed. For simplicity, we assume for the moment that $X$ is irreducible and reversible with respect to the uniform distribution on $I$. In this case, writing $\mathcal{L}$ for the generating matrix of $X$, we can conclude that the matrix $-\mathcal{L}$ is symmetric and admits an orthonormal basis of eigenvectors $v_1,v_2,\ldots,v_n$ corresponding to eigenvalues $0=\lambda_1<\lambda_2\leq\lambda_3\leq\ldots\leq\lambda_n$. Fixing a pair of initial states $(a,b)$ and letting $e_a$ (resp. $e_b$) be the vector, whose only non-zero component is the $a$-th one (resp. the $b$-th one) and equals to $1$, we have the decomposition
\eq
e_a-e_b=\sum_{l=2}^n \mu_l\,v_l.
\en 
Finally, we define the set
\eq
A(a,b):=\{0,\mu_2^2,\mu_2^2+\mu_3^2,\ldots,\mu_2^2+\mu_3^2+\ldots+\mu_n^2\}\subset[0,2]
\en
and a family of its neighborhoods
\eq
A^\delta(a,b):=[0,\delta\mu_2^2]\cup[(1-\delta)\mu_2^2,\mu_2^2+\delta\mu_3^2]\cup\ldots\cup[2-\delta\mu_n^2,2]
\en
for $\delta\in\Big(0,\frac{1}{2}\Big)$, and write $\|.\|_2$ for the $L^2$ norm with respect to the counting measure on $I$. 

\smallskip

\begin{thm}\label{thm5}
In the setting just described the following is true. For every $\delta\in\Big(0,\frac{1}{2}\Big)$, there is a constant $K(\delta)>0$ such that 
\eq
\int_0^\infty \mathbf{1}_{\big\{t\geq0:\;\|P^a_t-P^b_t\|_2^2\notin A^\delta(a,b)\big\}}\frac{1}{t}\,\mathrm{d}t\leq K(\delta)\,n
\en
for all pairs of initial states $a$, $b$. The constant $K(\delta)$ depends only on $\delta$, but not on $a$, $b$, $n$ or the particular Markov chain $X$.
\end{thm}

\medskip

The rest of the paper is organized as follows. In section 2.1, we prove a stronger version of Theorem \ref{thmTVgen} in the case that the invariant distribution $\nu$ is the uniform distribution on $I$. In order to do this, we show a novel entropy identity (see Lemma \ref{lemma2}) allowing us to relate the increase in the entropy of the one-dimensional distributions of the Markov chain to the convergence of the chain to its equilibrium. In section 2.2, we prove Theorem \ref{thmTVgen} by suitably adapting the entropy identity of Lemma \ref{lemma2} to the general setting. In section 3.1, we give a global (or, averaged) version of Theorem \ref{thm5} and present the proof of Theorem \ref{thm5}. Subsequently, we explain in section 3.2, how Theorem \ref{thm5} extends to general continuous time irreducible reversible finite-state Markov chains. Then, in section 4, we present a universal way of defining the ultrametric partition structure on the state space of a continuous time irreducible reversible finite-state Markov chain. Finally, in section 5, we show in the case that the chain is reversible with respect to the uniform distribution, how the entropy of the one-dimensional distributions of the chain can be used to obtain a control on the global convergence of the chain to its equilibrium.

\section{Estimates in total variation norm}

In this section we give a control on the convergence to equilibrium and the times of coupling with respect to the total variation norm by analyzing the change in the entropy of the Markov chain over time.  

\subsection{Markov chains reversible with respect to the uniform distribution}

The following theorem is a stronger version of Theorem \ref{thmTVgen} for the special case of Markov chains, which are reversible with respect to the uniform distribution. 

\begin{thm}\label{thm1}
Consider the setting of Theorem \ref{thmTVgen} and assume, in addition, that the Markov chain $X$ is reversible with respect to the uniform distribution on $I=\{1,2,\ldots,n\}$. Then: 
\begin{enumerate}[(a)]
\item There is a constant $C(\delta)>0$ depending only on $\delta$ (and not on $n$ or the particular Markov chain) such that for all initial states $a$ of the Markov chain:
\eq
\int_0^\infty \mathbf{1}_{\{t\geq0:\;\|P^a_{3t}-P^a_t\|_{TV}\geq\delta\}}\;\frac{1}{t}\;\mathrm{d}t \leq C(\delta)\log n.
\en
\item There is a constant $\tilde{C}(\delta)>0$ depending only on $\delta$ (and not on $n$ or the particular Markov chain) such that for all pairs $(a,b)$ of initial states of the Markov chain:
\eq
\int_0^\infty \mathbf{1}_{\{t\geq0:\;\|P^a_{t}-P^b_{t}\|_{TV}-\|P^a_{3t}-P^b_{3t}\|_{TV}\geq\delta\}}\;\frac{1}{t}\;\mathrm{d}t \leq \tilde{C}(\delta)\log n.
\en
\end{enumerate}
\end{thm}

\medskip

The proof relies on the following entropy identity.

\begin{lemma}\label{lemma2}
Let $X(t)$, $t\geq0$ be a Markov chain as in Theorem \ref{thm1} started in an initial state $a\in I$. Then, for all $t\geq0$:
\eq\label{ent_iden}
H(P^a_{2t})-H(P^a_t)=H(P^a_{t,2t}|P^a_{3t,2t}).
\en
Hereby, $H(.|.)$ stands for the relative entropy and $P^a_{u,s}$ stands for the law of the random vector $(X(u),X(s))$. 
In particular, the inequality
\eq\label{TVentropy}
\|P^a_t-P^a_{3t}\|_{TV}\leq\sqrt{2(H(P^a_{2t})-H(P^a_t))}
\en
holds for all $t\geq0$ and all initial states $a\in I$.
\end{lemma}

\noindent\textbf{Proof of Lemma \ref{lemma2}}. We start the proof with the following elementary computation, which only relies on the Markov property of $X$:
\begin{eqnarray*}
H(P^a_{2t})-H(P^a_t)&=& - \sum_{i\in I} P^a_{2t}(i)\,\log P^a_{2t}(i) + \sum_{i\in I} P^a_t(i)\,\log P^a_t(i) \\
&=& - \sum_{j\in I} P^a_t(j)\Big\{\Big(\sum_{i\in I} P^j_t(i)\,\log P^a_{2t}(i)\Big) - \log P^a_t(j) \Big\} \\
&=& - \sum_{j\in I} P^a_t(j)\sum_{i\in I} P^j_t(i)\Big\{\log P^a_{2t}(i)-\log P^a_t(j)\Big\} \\
&=& - \sum_{(i,j)\in I^2} P^a_t(j)\,P^j_t(i)\log\frac{P^a_{2t}(i)\,P^j_t(i)}{P^a_t(j)\,P^j_t(i)}.
\end{eqnarray*}
We now exploit the symmetry of the transition matrices of the Markov chain $X$ (which is due to the reversibility of the uniform distribution and the detailed balance condition) to deduce 
\eq
P^a_{2t}(i)\,P^j_t(i)=P^a_{2t}(i)\,P^i_t(j)
\en
for all $(i,j)\in I^2$. In addition, the Markov propery of $X$ yields
\begin{eqnarray}
&&P^a_{2t}(i)\,P^i_t(j)=\pp^a(X(2t)=i,X(3t)=j),\\
&&P^a_t(j)\,P^j_t(i)=\pp^a(X(t)=j,X(2t)=i)
\end{eqnarray}
for all $(i,j)\in I^2$. Putting the latter three observations together, we end up with the lemma. \ep

\bigskip

In the proof of Theorem \ref{thm1} we will need the following simple calculus lemma.

\begin{lemma}\label{lemma3}
Let $g:\,\rr\rightarrow[0,\infty)$ be a non-decreasing function, which satisfies
\eq
\lim_{u\rightarrow-\infty} g(u)=p,\quad\lim_{u\rightarrow\infty} g(u)=q
\en  
for some non-negative real constants $p\leq q$. Then, for every $r>0$ and $\epsilon>0$, one has the inequality
\eq
\int_{-\infty}^\infty \mathbf{1}_{\{g(u+r)-g(u)\geq\epsilon\}}\;\mathrm{d}u\leq \frac{r(q-p)}{\epsilon}\leq \frac{rq}{\epsilon}.
\en
\end{lemma}

\noindent\textbf{Proof of Lemma \ref{lemma3}.} It suffices to observe the elementary inequality
\eq
\mathbf{1}_{\{g(u+r)-g(u)\geq\epsilon\}}\leq\frac{g(u+r)-g(u)}{\epsilon},
\en
which leads to
\begin{eqnarray*}
\int_{-\infty}^\infty \mathbf{1}_{\{g(u+r)-g(u)\geq\epsilon\}}\;\mathrm{d}u
\leq\frac{1}{\epsilon}\lim_{K\rightarrow\infty}\Big(\int_{-K+r}^{K+r} g(u)\;\mathrm{d}u - \int_{-K}^{K} g(u)\;\mathrm{d}u \Big)
=\frac{r(q-p)}{\epsilon}
\end{eqnarray*}
and, hence, yields the lemma. \ep

\bigskip

We are now ready for the proof of Theorem \ref{thm1}. 
 
\bigskip 
 
\noindent\textbf{Proof of Theorem \ref{thm1}.} First, we note that part (b) of the theorem is a consequence of part (a) due to the inequalities
\eq
\|P^a_{t}-P^b_{t}\|_{TV}-\|P^a_{3t}-P^b_{3t}\|_{TV}\leq \|P^a_{t}-P^a_{3t}\|_{TV}+\|P^b_{t}-P^b_{3t}\|_{TV} 
\en 
and 
\eq\label{triangle}
\mathbf{1}_{\{\|P^a_{t}-P^a_{3t}\|_{TV}+\|P^b_{t}-P^b_{3t}\|_{TV}\geq\delta\}}
\leq\mathbf{1}_{\{\|P^a_{t}-P^a_{3t}\|_{TV}\geq\delta/2\}}
+\mathbf{1}_{\{\|P^a_{t}-P^a_{3t}\|_{TV}\geq\delta/2\}}. 
\en
We turn now to the proof of part (a). Due to the inequality \eqref{TVentropy}, it suffices to prove that for every $\overline{\delta}>0$ there is a constant $\overline{C}(\overline{\delta})>0$ depending only on $\overline{\delta}$ (and not on $a$, $n$ or the Markov chain $X$) such that
\eq
\int_0^\infty \mathbf{1}_{\{t\geq0:\;H(P^a_{2t})-H(P^a_t)\geq\overline{\delta}\}}\,\frac{1}{t}\;\mathrm{d}t\leq \overline{C}(\overline{\delta})\log n.
\en 
Introducing the function $g:\;\rr\rightarrow[0,\infty)$, $g(u)=H(P^a_{e^u})$, we can rewrite the latter inequality as
\eq
\int_{-\infty}^\infty \mathbf{1}_{\{u:\;g(u+\log 2)-g(u)\geq\overline{\delta}\}}\;\mathrm{d}u\leq \overline{C}(\overline{\delta})\log n.
\en
Noting that $\lim_{u\rightarrow\infty} g(u)=\lim_{t\rightarrow\infty} H(P^a_t)=\log n$ (since the uniform distribution is the unique stationary distribution of $X$), we see that the desired inequality holds with $\overline{C}(\overline{\delta})=\frac{\log 2}{\overline{\delta}}$ as a consequence of Lemma \ref{lemma3}. \ep 

\subsection{General reversible Markov chains}

In this subsection we consider a general continuous time irreducible reversible Markov chain $X$ on $I$ and will prove Theorem \ref{thmTVgen}. To start with, we recall the detailed balance condition:
\eq\label{det_bal}
\nu(a)\,\pp^a(X(t)=b)=\nu(b)\,\pp^b(X(t)=a),
\en
which holds for all times $t\geq0$ and all pairs of states $(a,b)\in I^2$. We now give the proof of Theorem \ref{thmTVgen}. 

\bigskip

\noindent\textbf{Proof of Theorem \ref{thmTVgen}.} The first assertion in part (b) of the theorem is a direct consequence of the inequality \eqref{triangle} (which clearly remains true in the more general setting of the present theorem) and the first assertion in part (a) of the theorem. Moreover, the second assertions in both parts of the theorem follow from the first assertions in the corresponding parts of the theorem and Markov's inequality. For these reasons, we only need to prove the first assertion in part (a) of the theorem. 

\bigskip

To this end, we fix an initial state $a\in I$ and note that the same computation as in the proof of Lemma \ref{lemma2} above yields:
\begin{eqnarray*}
H(P^a_{2t})-H(P^a_t)=\sum_{(i,j)\in I^2} P^a_t(j)\,P^j_t(i)\log\frac{P^a_t(j)\,P^j_t(i)}{P^a_{2t}(i)\,P^j_t(i)}.
\end{eqnarray*}
As before, we have by the Markov property
\eq
P^a_t(j)\,P^j_t(i)=\pp^a(X(t)=j,X(2t)=i).
\en
Moreover, the detailed balance condition \eqref{det_bal} gives
\eq
P^a_{2t}(i)\,P^j_t(i)=P^a_{2t}(i)\,P^i_t(j)\frac{\nu(i)}{\nu(j)}=\pp^a(X(2t)=i,\,X(3t)=j)\frac{\nu(i)}{\nu(j)}.
\en
Plugging this in, we get
\eq
H(P^a_{2t})-H(P^a_t)=H(P^a_{t,2t}|P^a_{3t,2t})+\sum_{(i,j)\in I^2} \pp^a(X(t)=j,X(2t)=i)\log\frac{\nu(i)}{\nu(j)},
\en 
where $P^a_{t,2t}$ and $P^a_{3t,2t}$ denote the laws of the random vectors $(X(t),X(2t))$ and $(X(3t),X(2t))$, conditioned on $X(0)=a$. In addition, writing $\log\frac{\nu(i)}{\nu(j)}=\log\nu(i)-\log\nu(j)$ and summing, we obtain
\begin{eqnarray*}
&&H(P^a_{2t})-H(P^a_t)-H(P^a_{t,2t}|P^a_{3t,2t})\\
&&=\sum_{i\in I} P^a(X(2t)=i)\log\nu(i)-\sum_{j\in I} P^a(X(t)=j)\log\nu(j).
\end{eqnarray*}  
Finally, integrating both sides of the latter equation with respect to $\nu$ and using the fact that $\nu$ is the invariant distribution of the Markov chain $X$, we end up with the \textit{averaged entropy identity}
\eq
\sum_{a\in I} \nu(a)\big(H(P^a_{2t})-H(P^a_t)\big)=\sum_{a\in I} \nu(a)\,H(P^a_{t,2t}|P^a_{3t,2t}). 
\en
In particular, this implies the inequality
\eq\label{gen_entr_TV}
\sum_{a\in I} \nu(a)\big(H(P^a_{2t})-H(P^a_t)\big)\geq \frac{1}{2}\sum_{a\in I} \nu(a)\,\|P^a_t-P^a_{3t}\|_{TV}^2.
\en

\medskip

On the other hand, the first inequality in part (a) of the theorem is equivalent to
\eq
\sum_{a\in I}\nu(a)\,\int_{-\infty}^\infty \mathbf{1}_{\{u\in\rr:\;\|P^a_{3e^u}-P^a_{e^u}\|_{TV}\geq\delta\}}\,\mathrm{d}u \leq C(\delta)\,H(\nu).
\en
This in turn would follow from $\mathbf{1}_{\{u\in\rr:\;\|P^a_{3e^u}-P^a_{e^u}\|_{TV}\geq\delta\}}\leq\frac{\|P^a_{3e^u}-P^a_{e^u}\|_{TV}^2}{\delta^2}$, if we can prove
\eq
\sum_{a\in I}\nu(a)\,\int_{-\infty}^\infty \frac{\|P^a_{3e^u}-P^a_{e^u}\|_{TV}^2}{\delta^2}\,\mathrm{d}u \leq C(\delta)\,H(\nu).
\en
However, due to the estimate \eqref{gen_entr_TV}, the left-hand side in the latter inequality is bounded above by
\begin{eqnarray*}
&&\sum_{a\in I}\nu(a)\,\int_{-\infty}^\infty \frac{2\big(H(P^a_{2e^u})-H(P^a_{e^u})\big)}{\delta^2}\,\mathrm{d}u\\
&=&\frac{2}{\delta^2}\sum_{a\in I}\nu(a)\lim_{K\rightarrow\infty}
\left(\int_{-K+\log 2}^{K+\log 2} H(P^a_{e^u})\,\mathrm{d}u-\int_{-K}^K H(P^a_{e^u})\,\mathrm{d}u\right)\\
&=&\frac{2\log 2\,H(\nu)}{\delta^2}.
\end{eqnarray*}
This finishes the proof. \ep

\section{Estimates in $L^2$ norm}

Throughout the first subsection of this section, we assume for the simplicity of notation that the continuous time Markov chain $X$ is irreducible and reversible with respect to the uniform distribution on the set $I=\{1,2,\ldots,n\}$. We first give a global version of Theorem \ref{thm5} in Theorem \ref{thm4} and then prove Theorem \ref{thm5} at the end of the first subsection. Subsequently, in the second subsection, we give the analogues of these results for a general continuous time irreducible reversible Markov chain on $I$.  

\subsection{Markov chains reversible with respect to the uniform distribution}

In the following theorem we show that, for most of the time on the scale-invariant clock, the square of the $L^2$ distance between the one-dimensional distributions of the Markov chain started in $a$ and the one-dimensional distributions of the Markov chain started in $b$, averaged over all pairs $(a,b)\in I^2$, stays close to the lattice
\eq
A_L:=\Big\{0,\frac{2}{n},\frac{4}{n},\ldots,\frac{2(n-1)}{n}\Big\}.
\en
This statement can be viewed as a global (or, averaged) version of Theorem \ref{thm5}. To make this statement precise, we write $A_L^\delta$ for the $\frac{2\delta}{n}$-neighborhood of $A_L$ in $\Big[0,\frac{2(n-1)}{n}\Big]$, where $\delta$ is a number in $\big(0,\frac{1}{2}\big)$, and can formulate the following result.    
\begin{thm}\label{thm4}
In the setting of Theorem \ref{thm5}, for all $0<\delta<\frac{1}{2}$, there exists a constant $K(\delta)>0$ such that the estimate
\eq
\int_0^\infty \mathbf{1}_{\big\{t\geq0:\;\frac{1}{n^2}\sum_{(a,b)\in I^2}\|P^a_t-P^b_t\|^2_2\notin A_L^\delta\big\}}
\frac{1}{t}\,\mathrm{d}t \leq K(\delta)\,n
\en
holds. Hereby, the constant $K(\delta)$ depends only on $\delta$, and not on $n$ or the particular Markov chain $X$. 
\end{thm}    

\bigskip

\noindent\textbf{Proof.} To start with, we recall the notation $\mathcal{L}$ for the generating matrix of the Markov chain $X$, so that, in particular, the transition matrix $P_t$ corresponding to a time $t\geq0$ is given by $e^{t\mathcal{L}}$. Since $X$ is irreducible and reversible with respect to the uniform distribution, the matrix $-\mathcal{L}$ is symmetric and non-negatively definite and has the eigenvalues $0=\lambda_1<\lambda_2\leq\lambda_3\leq\ldots\leq\lambda_n$. In particular, each of the matrices $P_t$, $t\geq0$ is symmetric, positively definite and has the eigenvalues 
\eq
1,\,e^{-\lambda_2 t},\,e^{-\lambda_3 t},\,\ldots,\,e^{-\lambda_nt}.
\en
Writing $\|.\|_2$ for the $L^2$ norm with respect to the counting measure on $I$ and $\langle.,.\rangle_2$ for the corresponding scalar product, we can make the following computation:
\begin{eqnarray*}
\sum_{(a,b)\in I^2}\|P^a_t-P^b_t\|_2^2&=&2n\sum_{a\in I} \|P^a_t\|_2^2 - 2\sum_{(a,b)\in I^2} \langle P^a_t,P^b_t\rangle_2\\
&=&2n\sum_{(a,c)\in I^2} \pp^a(X(t)=c)^2-2\Big\langle\sum_{a\in I} P^a_t,\sum_{b\in I} P^b_t\Big\rangle_2\\
&=&2n\big(e^{-2\lambda_2 t}+e^{-2\lambda_3 t}+\ldots+e^{-2\lambda_nt}\big),
\end{eqnarray*}
which is valid for all $t\geq0$. 

\bigskip

Next, we set 
\begin{eqnarray}
&&f(t):=e^{-2\lambda_2 t}+e^{-2\lambda_3 t}+\ldots+e^{-2\lambda_nt},\\
&&t_k:=\inf\{t\geq0:\;f(t)\leq k-\delta\},\quad k=1,\,2,\,\ldots,\,n-1.
\end{eqnarray}
The continuity of the function $f$ implies $f(t_k)=k-\delta$. In particular, it follows that
\eq
f(2t_k)\leq\max_{x_1+x_2+\ldots+x_{n-1}=k-\delta,\,0\leq x_i\leq 1} (x_1^2+x_2^2+\ldots+x_{n-1}^2).
\en 
Moreover, since the maximum of the convex function 
\[
(x_1,\ldots,x_{n-1})\mapsto(x_1^2+x_2^2+\ldots+x_{n-1}^2) 
\]
is taken over a convex set, it must be attained at a boundary point of that set. In other words, at the optimizing point it must hold $x_l\in\{0,1\}$ for at least one $1\leq l\leq n-1$. Eliminating the corresponding variable, we obtain a maximization problem of the same type and conclude that at least one another coordinate $x_{l'}$ has to belong to the set $\{0,1\}$. Proceeding with the same argument, we conclude that for each optimizing point $(x_1,x_2,\ldots,x_{n-1})$, there must be $(k-1)$ coordinates, which are equal to $1$, $(n-k-1)$ coordinates, which are equal to $0$, and one coordinate, which is equal to $1-\delta$. Thus, we have: 
\eq
f(2t_k)\leq f(t_k)-(1-\delta)+(1-\delta)^2=f(t_k)-\delta(1-\delta).
\en
Now, either $f(2t_k)\leq (k-1)+\delta$, or we can proceed with the same argument to conclude
\eq
f(4t_k)\leq f(2t_k)-\delta(1-\delta)\leq f(t_k)-2\delta(1-\delta).
\en
Proceeding further with the same argument, we end up with
\eq
f(2^R t_k)\leq f(t_k)-(1-2\delta)=(k-1)+\delta
\en
for $R=\Big\lceil\frac{1-2\delta}{\delta(1-\delta)}\Big\rceil$, where $\lceil.\rceil$ denotes the closest integer from above.  

\bigskip

Hence, setting
\eq
\tilde{t}_k=\inf\{t\geq0:\;f(t)\leq (k-1)+\delta\},\quad k=1,\,2,\,\ldots,\,n-1,
\en
we have the estimate
\eq
\int_0^\infty \mathbf{1}_{[t_k,\tilde{t}_k]}\frac{1}{t}\,\mathrm{d}t=\log\frac{\tilde{t}_k}{t_k}
\leq\log 2\,\Big\lceil\frac{1-2\delta}{\delta(1-\delta)}\Big\rceil=:K(\delta).
\en
Finally, using this and the identity 
\eq
\int_0^\infty \mathbf{1}_{\big\{t\geq0:\;\frac{1}{n^2}\sum_{(a,b)\in I^2}\|P^a_t-P^b_t\|^2_2\notin A_L^\delta\big\}}
\frac{1}{t}\,\mathrm{d}t=\sum_{k=1}^n \int_0^\infty \mathbf{1}_{[t_k,\tilde{t}_k]}\frac{1}{t}\,\mathrm{d}t, 
\en
we readily obtain the theorem. \ep

\medskip

\begin{rmk}
We note at this point that the main result of the article \cite{HP} implies that, for any vector $(1,\tilde{\lambda}_2,\tilde{\lambda}_3,\ldots,\tilde{\lambda}_n)$ with
\eq\label{mu_cond}
1>\tilde{\lambda}_2\geq\tilde{\lambda}_3\geq\ldots\geq\tilde{\lambda}_n>0, 
\en
there is a symmetric doubly stochastic matrix $S$ with eigenvalues $1,\tilde{\lambda}_2,\tilde{\lambda}_3,\ldots,\tilde{\lambda}_n$. 

\medskip

In particular, one can find a matrix $\mathcal{L}=S-\mathrm{Id}$ with the following two properties: 
\begin{enumerate}
\item[\textbf{(a)}] $\mathcal{L}$ generates a continuous time Markov chain, which is irreducible and reversible with respect to the uniform measure.
\item[\textbf{(b)}] The matrix $-\mathcal{L}$ has the eigenvalues
\eq
(0,\lambda_2,\lambda_3,\ldots,\lambda_n)=(0,1-\tilde{\lambda}_2,1-\tilde{\lambda}_3,\ldots,1-\tilde{\lambda}_n)
\en 
for a given sequence $\tilde{\lambda}_2,\tilde{\lambda}_3,\ldots,\tilde{\lambda}_n$ as in \eqref{mu_cond}.
\end{enumerate}

\medskip

This together with the proof of Theorem \ref{thm4} shows that the order $n$ of the upper bound in Theorem \ref{thm4} is optimal. As will become clear from the proofs below, the same is true for the upper bound of Theorem \ref{thm5}, and the counterparts of these results for general continuous time irreducible reversible Markov chains treated in section 3.2. 
\end{rmk}

\bigskip

We proceed with the proof of Theorem \ref{thm5}. 

\bigskip

\noindent\textbf{Proof of Theorem \ref{thm5}.} To start with, we introduce the following notations:
\begin{eqnarray*}
&&t_k(a,b):=\inf\Big\{t\geq0:\,\|P^a_t-P^b_t\|_2^2\leq\sum_{l=2}^{k-1} \mu_l^2 +(1-\delta)\mu_k^2\Big\},\; k=2,\ldots,n\\
&&\tilde{t}_k(a,b):=\inf\Big\{t\geq0:\,\|P^a_t-P^b_t\|_2^2\leq\sum_{l=2}^{k-1} \mu_l^2 +\delta\mu_k^2\Big\},\; k=2,\ldots,n
\end{eqnarray*}
and note that
\eq
\int_0^\infty \mathbf{1}_{\big\{t\geq0:\;\|P^a_t-P^b_t\|_2^2\notin A^\delta(a,b)\big\}}\frac{1}{t}\,\mathrm{d}t
=\sum_{k=2}^n \log\frac{\tilde{t}_k(a,b)}{t_k(a,b)}.
\en
From now on, we fix a $k=2,\,3,\,\ldots,\,n$ and will show that 
\eq\label{claim}
\log\frac{\tilde{t}_k(a,b)}{t_k(a,b)}\leq K(\delta) 
\en
for a suitable constant $K(\delta)>0$, which depends only on $\delta$ (but not on $a$, $b$, $k$ or $n$). To this end, we note that the identity
\eq
\|P^a_t-P^b_t\|_2^2=\sum_{l=2}^n \mu_l^2 e^{-2t\lambda_l},\quad t\geq0
\en
and the inequality $0<\lambda_2\leq\lambda_3\leq\ldots\leq\lambda_n$ imply the estimate
\eq
\|P^a_{2t_k(a,b)}-P^b_{2t_k(a,b)}\|_2^2 
\leq \max_{\sum_{l=2}^n \mu_l^2 x_l=\|P^a_{t_k(a,b)}-P^b_{t_k(a,b)}\|_2^2,\;1\geq x_2\geq x_3\geq\ldots\geq x_n\geq0} \sum_{l=2}^n \mu_l^2 x_l^2.
\en
Moreover, if we have $\mu_l^2>0$ for all $l=2,\,3,\,\ldots,\,n$, then the function
\[
(x_2,x_3,\ldots,x_n)\mapsto\sum_{l=2}^n \mu_l^2 x_l^2
\]
is stricly convex and must attain its maximum at a vertex point of the convex polyhedron
\[
\Big\{(x_2,x_3,\ldots,x_n):\;\sum_{l=2}^n \mu_l^2 x_l=\|P^a_{t_k(a,b)}-P^b_{t_k(a,b)}\|_2^2,\;1\geq x_2\geq x_3\geq\ldots\geq x_n\geq0\Big\}.
\]
If we have $\mu_l^2=0$ for some $l\in\{2,\,3,\,\ldots,\,n\}$, then we can elimininate the corresponding coordinate in the maximization problem and make the same conclusion for the reduced maximization problem. For this reason, we may assume without loss of generality that $\mu_l^2>0$ for all $l=2,\,3,\,\ldots,\,n$. Moreover, since the hyperplane $\sum_{l=2}^n \mu_l^2 x_l=\|P^a_{t_k(a,b)}-P^b_{t_k(a,b)}\|_2^2$ is $(n-2)$-dimensional, the vertices of the polyhedron above are given by points $1\geq x_2\geq x_3\geq\ldots\geq x_n\geq0$, for which $(n-2)$ of the inequalities
\[
1\geq x_2,\;x_2\geq x_3,\;\ldots,\;x_n\geq0
\] 
are in fact equalities. 

\bigskip

Thus, each optimizing point of the maximization problem above can be described as follows: There is a partition of $\{2,3,\ldots,n\}$ into three sets $I_1$, $I_2$, $I_3$ of the form $\{2,3,\ldots,l_1\}$, $\{l_1+1,l_1+2,\ldots,l_2\}$, $\{l_2+1,l_2+2,\ldots,n\}$, respectively, such that, for all $l\in I_1$ it holds $x_l=1$, for all $l\in I_2$ we have $x_l=\zeta$ for a suitable $\zeta\in[0,1]$, and for all $l\in I_3$ it holds $x_l=0$. Moreover, the identity
\eq
\|P^a_{t_k(a,b)}-P^b_{t_k(a,b)}\|_2^2=\sum_{l=2}^{k-1} \mu_l^2 + (1-\delta)\mu_k^2
\en    
shows that the value of $\zeta$ is given by 
\eq
\zeta=\frac{\big(\sum_{l=2}^{k-1} \mu_l^2 + (1-\delta)\mu_k^2\big) - \big(\sum_{l\in I_1} \mu_l^2\big)}{\sum_{l\in I_2}\mu_l^2},
\en
and that $I_1\subset\{2,3,\ldots,k-1\}$. To proceed, we introduce the set
\eq
\tilde{I}_2:=\big(\{2,3,\ldots,k-1\}\backslash I_1\big)\subset I_2
\en
and conclude 
\eq
\zeta=\frac{\big(\sum_{l\in \tilde{I}_2} \mu_l^2\big) +(1-\delta)\mu_k^2}{\sum_{l\in I_2}\mu_l^2},\quad \tilde{I}_2\subsetneq I_2.
\en
This allows us to make the following computation:
\begin{eqnarray*}
&&\|P^a_{2t_k(a,b)}-P^b_{2t_k(a,b)}\|_2^2
=\|P^a_{t_k(a,b)}-P^b_{t_k(a,b)}\|_2^2-\big(\sum_{l\in I_2} \mu_l^2 \zeta\big) + \big(\sum_{l\in I_2} \mu_l^2 \zeta^2\big)\\
&=&\|P^a_{t_k(a,b)}-P^b_{t_k(a,b)}\|_2^2
-\frac{\big(\sum_{l\in\tilde{I}_2}\mu_l^2 +(1-\delta)\mu_k^2\big)\big(\sum_{l\in I_2\backslash\tilde{I}_2} \mu_l^2 - (1-\delta)\mu_k^2\big)}
{\sum_{l\in I_2} \mu_l^2}. 
\end{eqnarray*}
Next, we note that the latter fraction is of the form $\frac{A\cdot B}{A+B}=\frac{1}{\frac{1}{A}+\frac{1}{B}}$, whereby: $A\geq(1-\delta)\mu_k^2$ and $B\geq \delta\mu_k^2$. Thus, 
\begin{eqnarray*}
\|P^a_{2t_k(a,b)}-P^b_{2t_k(a,b)}\|_2^2
&\leq& \|P^a_{t_k(a,b)}-P^b_{t_k(a,b)}\|_2^2 - \frac{1}{\frac{1}{(1-\delta)\mu_k^2}+\frac{1}{\delta\mu_k^2}}\\
&=& \|P^a_{t_k(a,b)}-P^b_{t_k(a,b)}\|_2^2 - \frac{\mu_k^2}{\frac{1}{1-\delta}+\frac{1}{\delta}}.
\end{eqnarray*}
Proceeding with the same argument, we conclude that
\eq
\|P^a_{2^R t_k(a,b)}-P^b_{2^R t_k(a,b)}\|_2^2 \leq \|P^a_{t_k(a,b)}-P^b_{t_k(a,b)}\|_2^2 - \mu_k^2 (1-2\delta) 
\en
for all natural numbers $R\geq\frac{1-2\delta}{\frac{1}{1-\delta}+\frac{1}{\delta}}$. In particular, we conclude that
\eq
\tilde{t}_k(a,b)\leq 2^{\big\lceil\frac{1-2\delta}{\frac{1}{1-\delta}+\frac{1}{\delta}}\big\rceil} t_k(a,b),
\en
where $\lceil.\rceil$ denotes the closest integer from above. This shows the claim \eqref{claim} with
\eq
K(\delta)=\log 2 \cdot \Big\lceil\frac{1-2\delta}{\frac{1}{1-\delta}+\frac{1}{\delta}}\Big\rceil
\en
and finishes the proof. \ep

\subsection{General reversible Markov chains}

We proceed with the analogues of Theorems \ref{thm5} and \ref{thm4} for a general continuous time irreducible reversible Markov chain $X$. To state the results, we introduce the following set of notations. We write $\nu$ for the invariant measure of $X$ as before, and let $D$ be the diagonal matrix, whose diagonal entries are given by $\nu(i)$, $i\in I$. Then, by the detailed balance condition \eqref{det_bal}, the matrix $D^{1/2}P_tD^{-1/2}$ is symmetric for all $t\geq0$. Moreover, since the matrices $D^{1/2}P_tD^{-1/2}$, $t\geq0$ commute, they have a joint orthonormal basis of eigenvectors $v_1,v_2,\ldots,v_n$ corresponding to sets eigenvalues
\eq
1>e^{-\lambda_2 t}\geq e^{-\lambda_3 t}\geq\ldots\geq e^{-\lambda_n t},\quad t\geq0,
\en
respectively (see chapter 3 of the book \cite{AF} for more details). In addition, for any fixed pair $(a,b)$ of initial states, we let 
\eq
D^{-1/2}(e_a-e_b)=\sum_{l=2}^n \tilde{\mu}_l v_l
\en
be the expansion of the vector $D^{-1/2}(e_a-e_b)$ in terms of the basis $v_1,v_2,\ldots,v_n$ (note that the vector $D^{-1/2}(e_a-e_b)$ is orthogonal to the eigenvector corresponding to the eigenvalue $1$ of the matrices $D^{1/2}P_tD^{-1/2}$, $t\geq0$). Finally, define the sets
\begin{eqnarray*}
\tilde{A}_L := \{0,1,2,\ldots,n-1\},\\
\tilde{A}(a,b) := \{0,\tilde{\mu}_2^2,\tilde{\mu}_2^2+\tilde{\mu}_3^2,\ldots,\tilde{\mu}_2^2+\tilde{\mu}_3^2+\ldots+\tilde{\mu}_n^2\},
\end{eqnarray*}
and their neighborhoods
\begin{eqnarray*}
\tilde{A}_L^\delta := [0,\delta]\cup[1-\delta,1+\delta]\cup[2-\delta,2+\delta]\cup\ldots\cup[n-1-\delta,n-1],\\ 
\tilde{A}^\delta(a,b) := [0,\delta\tilde{\mu}_2^2]\cup[(1-\delta)\tilde{\mu}_2^2,\tilde{\mu}_2^2+\delta\tilde{\mu}_3^2]\cup\ldots\cup
\Big[\sum_{l=2}^{n-1}\tilde{\mu}_l^2+(1-\delta)\tilde{\mu}_n^2,\sum_{l=2}^n \tilde{\mu}_l^2\Big],
\end{eqnarray*}
$0<\delta<\frac{1}{2}$. With these notations, the analogues of Theorems \ref{thm5} and \ref{thm4} read as follows.

\begin{thm}\label{thmL2gen}
Let $\|.\|_{L^2(\nu^{-1})}$ and $\langle.,.\rangle_{L^2(\nu^{-1})}$ be the $L^2$ norm and scalar product with respect to the weights $\nu(i)^{-1}$, $i\in I$. Then, for all $0<\delta<\frac{1}{2}$, there is a constant $\tilde{K}(\delta)>0$ such that the inequalities
\eq\label{gen_global}
\int_0^\infty \mathbf{1}_{\big\{t\geq0:\;\sum_{(a,b)\in I^2}\nu(a)\,\nu(b)\,\|P^a_t-P^b_t\|^2_{L^2(\nu^{-1})}\notin\tilde{A}_L^\delta\big\}}
\frac{1}{t}\,\mathrm{d}t \leq \tilde{K}(\delta)\,n
\en
and 
\eq\label{gen_local}
\int_0^\infty \mathbf{1}_{\big\{t\geq0:\;\|P^a_t-P^b_t\|_{L^2(\nu^{-1})}^2\notin \tilde{A}^\delta(a,b)\big\}}\frac{1}{t}\,\mathrm{d}t
\leq \tilde{K}(\delta)\,n, \quad (a,b)\in I^2,
\en
hold true. The constant $\tilde{K}(\delta)$ depends only on $\delta$, and not on $n$ or the particular Markov chain $X$.
\end{thm}

\medskip

\noindent\textbf{Proof.} In order to prove \eqref{gen_global}, we use the fact that $\nu$ is the invariant distribution of the Markov chain $X$ to deduce the identities
\begin{eqnarray*}
&&\sum_{(a,b)\in I^2}\nu(a)\,\nu(b)\,\|P^a_t-P^b_t\|^2_{L^2(\nu^{-1})}\\
&&=2\sum_{a\in I} \nu(a)\,\|P_t^a\|^2_{L^2(\nu^{-1})}-2\big\langle \sum_{a\in I} \nu(a)\,P_t^a,\sum_{b\in I} \nu(b)\,P_t^b\big\rangle\\
&&=2\sum_{a\in I} \nu(a) \sum_{c\in I} P_t^a(c)^2\nu(c)^{-1}-2\langle\nu,\nu\rangle_{L^2(\nu^{-1})}\\
&&=2\sum_{(a,c)\in I^2} \nu(a)\,P_t^a(c)^2\,\nu(c)^{-1}-2,
\end{eqnarray*}  
which hold for all $t\geq0$. Moreover, the latter sum is given by the sum of squares of the entries of the matrix $D^{1/2}P_tD^{-1/2}$ and is, hence, equal to $1+\sum_{l=2}^n e^{-2\lambda_l t}$. Thus,
\eq
\sum_{(a,b)\in I^2}\nu(a)\,\nu(b)\,\|P^a_t-P^b_t\|^2_{L^2(\nu^{-1})}=\sum_{l=2}^n e^{-2\lambda_l t}\,\quad t\geq0.
\en
From this point on, one can proceed as in the proof of Theorem \ref{thm4} to show \eqref{gen_global}.

\bigskip

Now, we turn to the proof of \eqref{gen_local}. To this end, we note that the detailed balance condition \eqref{det_bal} implies $DP_t=P_t^T D$, $t\geq0$, where the superscript $T$ stands for the transpose of a matrix. This allows us to make the computation
\begin{eqnarray*}
P_t^a-P_t^b=((e_a-e_b)^T P_t)^T=P_t^T (e_a-e_b)=D^{1/2}\big(D^{1/2}P_tD^{-1/2}\big)D^{-1/2}(e_a-e_b)\\
= D^{1/2}\sum_{l=2}^n \tilde{\mu}_l\, e^{-\lambda_l t}\, v_l
= \sum_{l=2}^n \tilde{\mu}_l\, e^{-\lambda_l t}\, D^{1/2} v_l
\end{eqnarray*}
for all $t\geq0$. Next, we observe that the vectors $D^{1/2}v_1,D^{1/2}v_2,\ldots,D^{1/2}v_n$ form an orthonormal basis with respect to the scalar product $\langle.,.\rangle_{L^2(\nu^{-1})}$, since the vectors $v_1,v_2,\ldots,v_n$ form an orthonormal basis with respect to the standard Euclidean scalar product. Hence,
\eq
\|P_t^a-P_t^b\|^2_{L^2(\nu^{-1})}=\sum_{l=2}^n \tilde{\mu}_l^2\,e^{-2\lambda_l t},\quad t\geq0.
\en
From this point on, one only needs to follow the arguments in the proof of Theorem \ref{thm5} to end up with \eqref{gen_local}. \ep

\begin{rmk}
It is worth noting that the estimates of Theorems \ref{thm5}, \ref{thm4} and \ref{thmL2gen} hold for 
\[
\|P_t^a-\nu\|_2^2,\;\;\frac{2}{n}\sum_{a\in I} \|P^a_t-\nu\|_2^2,\;\;2\sum_{a\in I} \nu(a)\|P^a_t-\nu\|_{L^2(\nu^{-1})}^2\;\;\text{and}\;\;
\|P_t^a-\nu\|_{L^2(\nu^{-1})}^2 
\]
in place of 
\begin{eqnarray*}
\|P_t^a-P_t^b\|_2^2,\;\;\frac{1}{n^2}\sum_{(a,b)\in I^2} \|P^a_t-P^b_t\|_2^2,\;\;\sum_{(a,b)\in I^2} \nu(a)\nu(b) \|P^a_t-P^b_t\|_{L^2(\nu^{-1})}^2 \\
\text{and}\;\;\|P^a_t-P^b_t\|_{L^2(\nu^{-1})}^2, 
\end{eqnarray*}
respectively. The same proofs apply, with the only difference being that one needs to expand the vectors $(e_a-\nu)$ and $D^{-1/2}(e_a-\nu)$ in terms of an orthonormal basis of eigenvectors of the matrices $P_t$, $t\geq0$ and $D^{1/2}P_tD^{-1/2}$, $t\geq0$, respectively.
\end{rmk}

\section{A universal approach to the ultrametric structure}

In this section we provide a univeral way of defining the ultrametric partition structure on the state space $I=\{1,2,\ldots,n\}$ of a continuous time irreducible Markov chain $X$, which is reversible with respect to its invariant distibution $\nu$. Typical examples of such chains are encountered in statistical physics, where often the transition rate for a pair $(a,b)$ of neighboring states is proportional to $e^{-\beta(E(b)-E(a))_+}$ with $E$ being an energy functional (see the references given in the introduction, as well as the references therein). For large values of $\beta$, the energy landscape naturally provides a partition of the state space into states of different types, which are separated by potential wells (see Figure 1 for a schematic diagram). 

\begin{figure}[h]
\includegraphics[width=10cm, height=7cm]{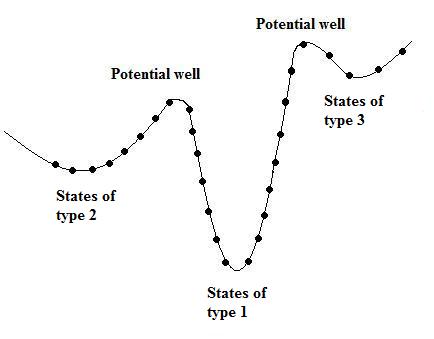}
\caption{A schematic diagram of an energy landscape}
\end{figure}

Here, we will give a universal way of defining the partition structure without making use of the explicit knowledge of the transition rates. Thereby, each of the partitions will correspond to a time scale on which convergence to equilibrium occurs for the Markov chain in consideration. For this purpose, we let $\|.\|$ be any norm on the space of finite measures on the set $I$, which is normalized in such a way that $\|\pi_1-\pi_2\|\leq 1$ for any two probability measures $\pi_1$, $\pi_2$ on $I$. Moreover, we assume that the function $t\mapsto\|\pi_1 P_t-\pi_2 P_t\|$ is strictly decreasing on $[0,\infty)$ and tends to zero in the limit $t\rightarrow\infty$ for all probability measures $\pi_1\neq\pi_2$ on $I$ (hereby, the products $\pi_1 P_t$, $\pi_2 P_t$ should be understood in the sense of multiplication of a probability measure by a stochastic kernel). Examples of such norms are the appropriately normalized total variation and $L^2$ norms discussed above.

\bigskip

Now, we fix an $0<\epsilon<1$ and will recursively define equivalence relations $\sim_1, \sim_2, \ldots$ on $I$, which will induce the desired sequence of nested partitions. To define $\sim_1$, we set
\begin{eqnarray}
&&t^*_1=\min_{a\neq b}\;\inf\{t\geq0:\;\|P_t^a-\nu\|+\|P_t^b-\nu\|\leq\epsilon\},\\
&&h^{(1)}_a(t)=\|P_{t^*_1}^a-\nu\|^k\cdot\|P_s^a-\nu\|,\quad\text{for}\;t=kt^*_1+s,\;0\leq s<t^*_1,\;a\in I. 
\end{eqnarray}   
Then, we let $a\sim_1 b$ iff either $a=b$, or
\eq
\limsup_{t\rightarrow\infty}\frac{t^*_1}{t}\log\big[h^{(1)}_a(t)+h^{(1)}_b(t)\big]\leq \log (2\epsilon).
\en
Now, to define $\sim_2$, we set
\begin{eqnarray}
&&t^*_2=\min_{a\not\sim_1 b}\;\inf\{t\geq0:\;\|P_t^a-\nu\|+\|P_t^b-\nu\|\leq\epsilon\},\\
&&h^{(2)}_a(t)=\|P_{t^*_2}^a-\nu\|^k\cdot\|P_s^a-\nu\|,\quad\text{for}\;t=kt^*_2+s,\;0\leq s<t^*_2,\;a\in I.
\end{eqnarray}
Then, we let $a\sim_2 b$ iff either $a\sim_1 b$, or
\eq
\limsup_{t\rightarrow\infty}\frac{t^*_2}{t}\log\big[h^{(2)}_a(t)+h^{(2)}_b(t)\big]\leq \log (2\epsilon).
\en
The equivalence relations $\sim_3, \sim_4, \ldots$ are now defined analogously.

\bigskip

The intuition behind the definitions above can be explained as follows. For each $l\in\nn$, the time $t^*_l$ is defined as the first time, at which there is a pair of states $(a,b)$, which have not already been declared to be equivalent with respect to $\sim_{l-1}$ and for which both the distance of $P^a_t$ and the distance of $P^b_t$ from the equilibrium distribution $\nu$ is small. For such a pair $(a,b)$ the identity $a\sim_l b$ is due to the following computation:
\eq
\limsup_{k\rightarrow\infty}\frac{t^*_l}{k t^*_l}\log\big[h^{(l)}_a(kt^*_l)+h^{(l)}_b(kt^*_l)\big]
\leq\limsup_{k\rightarrow\infty}\frac{t^*_l}{k t^*_l}\log (2\epsilon^k)=\log\epsilon. 
\en
Increasing the right-hand side of the inequality defining $\sim_l$ to $\log (2\epsilon)$ allows us to find the pairs of states $(c,d)$, for which the distributions $P^c_t$, $P^d_t$ approach the equilibrium distribution $\nu$ on approximately the same time scale as $P^a_t$, $P^b_t$. The functions $h^{(l)}_c$, $h^{(l)}_d$ are hereby, in a suitable sense, our best guess for the functions $t\mapsto\|P^c_t-\nu\|$, $t\mapsto\|P^d_t-\nu\|$, if we only observe the latter on the time interval $[0,t^*_l]$. The following proposition summarizes our findings.

\begin{prop}
The relations $\sim_1, \sim_2, \ldots$ defined above are equivalence relations and define a sequence of nested partitions of the state space $I=\{1,2,\ldots,n\}$. Moreover, it holds $a\sim_l b$ for any pair $(a,b)$ which achieves the minimum in
\eq\label{whatistl}
t^*_l=\min_{a\not\sim_{l-1} b} \inf\{t\geq0:\;\|P_t^a-\nu\|+\|P_t^b-\nu\|\leq\epsilon\}
\en
and we have $a\sim_{n-1} b$ for any pair $(a,b)\in I^2$. 
\end{prop} 

\medskip

\noindent\textbf{Proof.} Fix an $l\in\nn$. To show that $\sim_l$ is an equivalence relation, we only need to prove the transitivity of $\sim_l$. To this end, we observe that the inequality
\eq
[h^{(l)}_a(t)+h^{(l)}_c(t)]\leq [h^{(l)}_a(t)+h^{(l)}_b(t)] + [h^{(l)}_b(t)+h^{(l)}_c(t)],\quad t\geq0
\en
together with Lemma 1.2.15 in Chapter 1 of \cite{DZ} yield
\begin{eqnarray*}
&&\limsup_{t\rightarrow\infty}\frac{1}{t}\log\big[h^{(l)}_a(t)+h^{(l)}_c(t)\big]\\
&&\leq \max\Big(\limsup_{t\rightarrow\infty}\frac{1}{t}\log\big[h^{(l)}_a(t)+h^{(l)}_b(t)\big],\;
\limsup_{t\rightarrow\infty}\frac{1}{t}\log\big[h^{(l)}_b(t)+h^{(l)}_c(t)\big]\Big)
\end{eqnarray*}
for all $(a,b,c)\in I^3$. Hence, the relations $a\sim_l b$ and $b\sim_l c$ imply together $a\sim_l c$. Moreover, since $a\sim_{l-1}b$ implies $a\sim_l b$ by definition, and $a\sim_l b$ holds for each pair $(a,b)\in I^2$, which achieves the minimum in \eqref{whatistl} (see the paragraph preceeding the proposition), the number of equivalence classes under $\sim_l$ is at most $n-l$. This shows $a\sim_{n-1} b$ for all pairs $(a,b)\in I^2$. \ep

\section{Bounds on the global convergence to equilibrium through the entropy}

We have seen in section 2 that one can obtain a control on the convergence to equilibrium and the times of coupling by analyzing the entropy that is accumulated by the Markov chain over time. In this section, we pursue this idea further and give estimates on the approach to equilibrium on subsets of macroscopic size for continuous time irreducible Markov chains which are \textit{reversible with respect to the uniform distribution}. To this end, for each $0<\kappa<1$ and $t\geq0$, we introduce the set
\eq
EQ(t)=\Big\{a\in I:\;P_t(a)\in\Big(\frac{1-\kappa}{n},\frac{1+\kappa}{n}\Big)\Big\},
\en
where, with a slight abuse of notation, we wrote $P_t$ for the law of the random variable $X(t)$. For each $t\geq0$, the set $EQ(t)\subset I$ should be viewed as the part of the state space on which the probability measure $P_t$ is close to the equilibrium distribution of the Markov chain $X$. We are interested in lower bounds on the size $|EQ(t)|$ of such sets.

\begin{thm}\label{thm2}
Fix real numbers $0<\kappa<1$ and $0<\alpha<\frac{1}{2}$, set $\tilde{\alpha}=\frac{1}{2}-\alpha$, and, on the interval $[0,\tilde{\alpha}]$, define the function
\begin{eqnarray*}
F(\alpha_1)&=&-\alpha_1(1 -\kappa)\log(1 - \kappa) - (\tilde{\alpha}-\alpha_1)(1 + \kappa)\log(1 + \kappa) \\
&&- (1 - (1 - \kappa)\alpha_1 - (1 + \kappa)(\tilde{\alpha}-\alpha_1))
\log\left(\frac{1-\tilde{\alpha}-\tilde{\alpha}\kappa+2\kappa\alpha_1}{1-\tilde{\alpha}}\right)
\end{eqnarray*}
taking non-positive values. Then, the entropy estimate 
\begin{equation}\label{entr_cond}
H(P_t)>\log n + \max_{0\leq\alpha_1\leq\tilde{\alpha}} F(\alpha_1)
\end{equation}
implies the lower bound
\eq\label{E_LBD}
|EQ(t)|\geq\alpha n.
\en
Hereby, depending on the values of $\kappa$ and $\alpha$, the maximum in \eqref{entr_cond} is attained at $0$, $\tilde{\alpha}$ or
\[
\alpha_1^*:=\frac{(1-\kappa)^{-1/(2\kappa)}
\Big((1-\tilde{\alpha})\sqrt{1-\kappa}(1+\kappa)^{(1+\kappa)/(2\kappa)}-e(1-\kappa)^{1/(2\kappa)}(1-\tilde{\alpha}-\tilde{\alpha}\kappa)\Big)}
{2e\kappa}. 
\]
\end{thm}

\medskip

\noindent\textbf{Proof.} We fix numbers $\kappa$ and $\alpha$ as in the statement of the theorem and suppose that the inequality \eqref{E_LBD} does not hold. We will show that this implies that the entropy bound \eqref{entr_cond} cannot hold. To start with, we introduce the notation $p_a:=P_t(a)$, $a\in I$, and make the decomposition
\eq
H(P_t)=-\sum_{a\in EQ(t)} p_a\log p_a - \sum_{b\notin EQ(t)} p_b\log p_b. 
\en 

\medskip
 
For a given value of $\rho:=\sum_{a\in EQ(t)} p_a \in [0,1]$, the maximum of the function $-\sum_{b\notin EQ(t)} p_b\log p_b$ is attained on the interior boundary of the set 
\eq\label{opt_set}
\Big\{\sum_{b\notin EQ(t)} p_b = 1-\rho:\;p_b\notin\Big(\frac{1-\kappa}{n},\frac{1+\kappa}{n}\Big)\Big\}.  
\en
Indeed, this is a consequence of the fact that the function 
\eq\label{ent_funct}
(p_b:\;b\notin EQ(t))\mapsto - \sum_{b\notin EQ(t)} p_b\log p_b
\en
is concave and attains its maximum over the convex set $\{\sum_{b\notin EQ(t)} p_b = 1-\rho\}$ at the point $(\frac{1-\rho}{n-|EQ(t)|},\frac{1-\rho}{n-|EQ(t)|},\ldots,\frac{1-\rho}{n-|EQ(t)|})$, which is not an element of the set in \eqref{opt_set}. The latter statement follows from the inequalities
\eq
\frac{1-\rho}{n-|EQ(t)|}\geq\frac{1-\frac{1+\kappa}{n}|EQ(t)|}{n-|EQ(t)|}>\frac{1-\kappa}{n}
\en
and
\eq
\frac{1-\rho}{n-|EQ(t)|}\leq\frac{1-\frac{1-\kappa}{n}|EQ(t)|}{n-|EQ(t)|}<\frac{1+\kappa}{n}
\en
with the respective second inequalities in the latter two displays being consequences of $|EQ(t)|<\frac{n}{2}$. 

\bigskip

From the preceeding argument we conclude that at least one of the coordinates of a point in the set \eqref{opt_set}, which maximizes the function in \eqref{ent_funct}, has be equal to $\frac{1-\kappa}{n}$ or $\frac{1+\kappa}{n}$. Eliminating this coordinate and proceeding with the same argument, we deduce that at least $\frac{n}{2}-|EQ(t)|$ coordinates of an optimizing point have to be equal to $\frac{1-\kappa}{n}$ or $\frac{1+\kappa}{n}$. Now, eliminating all coordinates, which belong to the set $\Big\{\frac{1-\kappa}{n},\frac{1+\kappa}{n}\Big\}$, we deduce the following: If the inequality \eqref{E_LBD} fails, then the entropy $H(P_t)$ cannot exceed the entropy of a probability measure on a set of $n$ elements, for which at least $\frac{n}{2}-|EQ(t)|$ of its weights belong to the set $\Big\{\frac{1-\kappa}{n},\frac{1+\kappa}{n}\Big\}$ and the rest of its weights is equal. In other words, denoting the proportion of weights, which are equal to $\frac{1-\kappa}{n}$, by $\alpha_1$ and the proportion of weights, which are equal to $\frac{1+\kappa}{n}$, by $\alpha_2$, we have: $H(P_t)\leq\max_{\alpha_1,\alpha_2}\overline{F}(\alpha_1,\alpha_2)$ with
\begin{eqnarray*}
\overline{F}(\alpha_1,\alpha_2)&=&\Big[-\alpha_1(1-\kappa)\log\frac{1-\kappa}{n}-\alpha_2(1+\kappa)\log\frac{1+\kappa}{n}\\
&&-(1-(1-\kappa)\alpha_1-(1+\kappa)\alpha_2)\log\frac{1-(1-\kappa)\alpha_1-(1+\kappa)\alpha_2}{(1-\alpha_1-\alpha_2)n}\Big].
\end{eqnarray*}  
Hereby, the maximum is taken under the constraints $\frac{1}{2}-\alpha\leq\alpha_1+\alpha_2\leq1$, $\alpha_1\geq0$, $\alpha_2\geq0$, $1-(1-\kappa)\alpha_1-(1+\kappa)\alpha_2\geq 0$.

\bigskip

Next, we note that $\overline{F}$ can be written as $(\log n)+\tilde{F}$, where $\tilde{F}$ is given by
\begin{eqnarray*}
\tilde{F}(\alpha_1,\alpha_2)&=&\Big[-\alpha_1(1-\kappa)\log(1-\kappa)-\alpha_2(1+\kappa)\log(1+\kappa)\\
&&-(1-(1-\kappa)\alpha_1-(1+\kappa)\alpha_2)\log\frac{1-(1-\kappa)\alpha_1-(1+\kappa)\alpha_2}{1-\alpha_1-\alpha_2}\Big].
\end{eqnarray*}  
Hence, $H(P_t)\leq(\log n)+\max_{\alpha_1,\alpha_2}\tilde{F}(\alpha_1,\alpha_2)$, where the maximum is taken over the region described at the end of the previous paragraph. Now, a straightforward computation of the Hessian of $\tilde{F}$ together with the constraint $1-(1-\kappa)\alpha_1-(1+\kappa)\alpha_2\geq 0$ show that the function $\tilde{F}$ is concave throughout the region over which its maximum is taken. In addition, the maximum of $\tilde{F}$ over the region determined by the constraints $\alpha_1\geq0$, $\alpha_2\geq0$, $\alpha_1+\alpha_2\leq1$, $1-(1-\kappa)\alpha_1-(1+\kappa)\alpha_2\geq 0$ is attained at the point $(0,0)$ and is equal to $0$, since it corresponds to the highest value of the entropy that a probability measure on a set of $n$ elements can take (namely, $\log n$). Thus, the maximum of $\tilde{F}$ over the region of interest is attained on the segment given by the constraints $\alpha_1\geq0$, $\alpha_2\geq0$, $\alpha_1+\alpha_2=\frac{1}{2}-\alpha$. Plugging in $\frac{1}{2}-\alpha-\alpha_1$ instead of $\alpha_2$, and recalling the notation $\tilde{\alpha}=\frac{1}{2}-\alpha$, we end up with $H(P_t)\leq(\log n)+\max_{0\leq\alpha_1\leq\tilde{\alpha}} F(\alpha_1)$. This is the desired contradiction to \eqref{entr_cond}. 

\bigskip

We also observe that the function $F$ must be non-positive throughout $[0,\tilde{\alpha}]$, since the entropy of a probability measure on a set of $n$ elements cannot exceed the value $\log n$. Moreover, since the function $\tilde{F}$ is concave, the function $F$ is also concave. Furthermore, a straightforward computation shows that, depending on the values of $\kappa$ and $\alpha$, either the derivative of the function $F$ has no zeros on the interval $[0,\tilde{\alpha}]$, in which case $F$ attains its maximum at one of the boundary points, or the only zero of the derivative of $F$ on the interval $[0,\tilde{\alpha}]$ is given by $\alpha_1^*$ (defined in the statement of the theorem), in which case $F$ attains its maximum at $\alpha_1^*$. This finishes the proof. \ep  

\section{Acknowledgement}

The author would like to thank David J. Aldous for his comments throughout the preparation of this work. He is also grateful to Anton Bovier for his remarks on an early version of this manuscript.

\bigskip\bigskip

\bibliographystyle{alpha}

\bigskip\bigskip\bigskip\noindent

\end{document}